\documentclass[11pt,twoside,reqno]{amsart}
\usepackage{mathrsfs}
\usepackage{amsmath}
\usepackage{amssymb}
\usepackage[a4paper]{geometry}
\theoremstyle{plain}
\newtheorem{teo}{Theorem}
\newtheorem{pro}{Proposition}
\newtheorem{lem}{Lemma}
\newtheorem{cor}{Corollary}

\theoremstyle{definition}
\newtheorem{defi}{Definition}
\theoremstyle{remark}

\newtheorem{nota}{Remark}

\newcommand{\tr}{\mathop{\rm tr}}

\newcommand{\re}{\mathbb{R}}

\newcommand{\g}{\Gamma_m}
\newcommand{\f}{{}_0F_1}
\newcommand{\fc}{{}_1\mathscr{F}_1}

\newcommand{\df}{\tilde{H_m^+}}
\newcommand{\h}{{}_0\mathscr{F}_1}
\date{\today}
\title{Laguerre Process and Generalized Hartman-Watson Law } 
\begin{document}
\maketitle
\centerline{N. DEMNI$^{\star}$}
 
 ${}^{\star}$ {\footnotesize Laboratoire de Probabilit\'es et Mod\`eles Al\'eatoires, Universit\'e de Paris VI, 4 Place Jussieu, Case 188, F-75252 Paris Cedex 05,
 e-mail : demni@ccr.jussieu.fr 
 \\ {\bf Key Words}:  Generalized Hartman-Watson law, Gross-Richards formula, Laguerre process, Special functions of matrix argument.} 
 
\begin{abstract} : In this paper, we study  complex Wishart processes or the so-called \emph{Laguerre processes} $(X_t)_{t \geq 0}$. We are interested in the behaviour of the eigenvalue process, we derive some useful stochastic  differential equations and  compute both the infinitesimal generator and the semi-group. We also give  absolute-continuity relations between different indices. Finally, we compute the density function of the so-called {\it generalized Hartman-Watson} law as well as the law of $T_0 := \inf\{t, \det(X_t) = 0\}$ when the size of the matrix is $2$.\end{abstract}

\section{Introduction}
The Real Wishart process is a symmetric matrix-valued process which was introduced by M.F.Bru (1989) as follows : Let $B_t = (B_{ij}(t))_{i,j}$ be a $n \times m$ Brownian matrix and define 
$X_t = B_t^T B_t$. The process $(X_t)_{t \geq 0}$ satisfies the following stochastic differential equation (SDE)
\begin{align*}
dX_t & = B_t^TdB_t + dB_t^T B_t  + n{\it I}_m dt 
\\& = \sqrt{X_t}dN_t + dN_t^T \sqrt{X_t} + n{\it I}_m dt , \qquad X_0 = B_0^TB_0
\end{align*} where ${\it I}_m$ denotes the unit matrix, the superscript ${}^T$ stands for the transpose, $\sqrt{X_t}$ is the matrix square root of the positive matrix $X_t$ and $(N_t)_{t \geq 0}$ is a $m \times m$ Brownian matrix. This process is called the Wishart process of dimension $n$, of size $m$, starting from $X_0$ and is denoted $W(n,m,X_0)$. Then, the $W(\delta,m,X_0)$ where $\delta$ runs over the Gindikin ensemble $\{1,\dots, m-1\} \cup]m-1,\infty[$ is defined as the unique solution of the latter SDE with $\delta$ instead of $n$. Thus, it can be viewed as an extension of the squared Bessel process to higher dimension. In this way,  Donati et al. (2004) tried to extend some well known properties of the squared Bessel processes to the matrix case and derived expressions such as the Laplace transform and the tail distribution of some random variables, in which many multivariate special functions of symmetric matrix argument appear, such as Gamma, modified Bessel and hypergeometric functions (Muirhead 1982). However, the latter is quite complicated to deal with and to our knowledge, there are no more precise results on the law of these variables.  
Nevertheless, in the complex case, hypergeometric functions of Hermitian matrix argument can be expressed as a determinant of a matrix whose entries are one-dimensional hypergeometric functions. In fact, Gross and Richards (1991) established the following result :
\begin{equation*}
{}_pF_q(a_1, \hdots, a_p, b_1, \hdots, b_q; X) = \frac{\det(x_i^{n-j}{}_p\mathscr{F}_q(a_1 - j +1,\hdots, a_p - j +1, \hdots, b_q - j + 1; x_i))}{V(X)}
\end{equation*} where $X$ is a $m \times m$ Hermitian matrix, $(x_i)$ are its eigenvalues, ${}_p\mathscr{F}_q$ denotes the standard hypergeometric functions with scalar argument, 
$V(X) = \prod_{i < j}(x_i - x_j)$ is the Vandermonde determinant and ${}_pF_q$ is the hypergeometric function with Hermitian matrix argument defined by :  
\begin{equation*} 
{}_pF_q(a_1, \hdots, a_p, b_1, \hdots, b_q; X) = \sum_{k \geq 0}\sum_{\tau}\frac{(a_1)_{\tau}\cdots(a_p)_{\tau}}{(b_1)_{\tau}\cdots(b_q)_{\tau}} \frac{C_{\tau}(X)}{k!}
\end{equation*}where $\tau = (\tau_1, \hdots, \tau_m)$ is a partition of length $\leq m$ and of weight $k$ (i.e.  $\tau_1 \geq \tau_2 \geq \dots \geq \tau_m, \sum_i\tau_i = k$), $(a)_{\tau}$ is the generalized Pochammer symbol and $C_{\tau}$ is the so-called zonal polynomial. We refer to Macdonald (1995) for further details and Lassalle (1991), (1991a) for analogous expressions for  multivariate orthogonal polynomials. The determinantal representation above is due to the fact that the zonal polynomial is identified with the (normalized) \emph{Schur functions} defined by :  
\begin{equation*}   
s_{\tau}(x_1, \hdots, x_m) = \frac{\det(x_i^{\tau_j + m - j})}{\det(x_i^{m-j})}
\end{equation*}
  
Consequently, one can use integral representations as well as other properties of standard hypergeometric functions to get, at least when $m=2$, some results that are till now
unknown in the Wishart case.    
The rest of this paper consists of seven sections, which are respectively devoted to the following topics: in section $2$, we introduce the Laguerre process of integer dimension. In section $3$, we study the behaviour of the eigenvalue process. Then, in section $4$, we define the Laguerre process of positive real dimension. Section $5$ is devoted to the absolute-continuity relations, from which we deduce the Laplace transform of the so-called {\it generalized Hartman-Watson} law as well as the tail distribution of $T_0$, the first hitting time of $0$. 
In section $6$, we focus on the case $m=2$ for which we invert this Laplace transform, and finally, in section $7$, we compute the density of $S_0 : = 1/(2T_0)$. 


\section{Laguerre Process of integer index}
Let B be a $n \times m$ complex Brownian matrix starting from $B_0$ , i.e. $B = (B_{ij})$ where the entries $B_{ij}$ are independent complex Brownian motions, so we can write $B = B^1 + iB^2$ where $B^1, B^2$ are two independent real Brownian matrices. We are interested in the matrix-valued  process $X_t := B_t^{\star}B_t$ which satisfies the following SDE :
\begin{equation}
 dX_t = dB_t^{\star}B_t + B_t^{\star}dB_t + 2n{\it I}_mdt \label{Henrike}\end{equation} 
\begin{defi}  $(X_t)_{t \geq 0}$ is called the Laguerre process of size $m$, of  dimension $n$ and starting from $X_0 = B_0^{\star}B_0$, and will be denoted by $L(n,m,X_0)$.\end{defi}
\begin{nota}
For $m = 1$, $(X_t)_{t \geq 0}$ is a squared Bessel process of dimension $2n$, denoted by $BESQ(2n,X_0)$.\end{nota}
\begin{nota} Set $X_t = (X_{ij}(t))_{i, j}$. One can easily check that 
\begin{equation*} dX_{ii}(t) = 2\sqrt{X_{ii}(t)} d\gamma_i(t) + 2ndt \qquad 1 \leq i \leq m, \end{equation*}   
 where $(\gamma_i)_{1 \leq i \leq m}$ are independent Brownian motions, thus, $X_{ii}$ is a $BESQ (2n,(X_0)_{ii})$.
 \end{nota}
\begin{nota}
The equation above implies that :
\begin{equation} d(\mathop{\rm tr}(X_t)) = 2 \sqrt{\tr(X_t)} d\beta_t  + 2nm dt
\end{equation} where $\beta$ is a Brownian motion. Consequently, $(\tr(X_t))_{t \geq 0}$ is a $BESQ(2nm,\tr(X_0))$ of dimension $2nm$ starting from $\tr(X_0)$. One can also deduce from (\ref{Henrike}) that for every $i, j, k, l \in \{1, \cdots, m\}$: 
\begin{equation*} 
\langle dX_{ij}, dX_{kl}\rangle_t = 2(X_{il} \delta_{kj} + X_{kj} \delta_{il})dt,   \end{equation*} 
which differs from equation (I-1-5) derived by Bru (1989) since for a complex Brownian motion $\gamma$, we have $d\langle \gamma, \gamma \rangle_t = 0$ and $\langle\gamma,\overline{\gamma} \rangle_t = 2t$ .
\end{nota}
\subsection{Infinitesimal generator}
Let $H_m,\df$ be respectively the space of $m \times m$ Hermitian matrices and the space of $m \times m$ positive definite Hermitian matrices. On the space of Hermitian matrix-argument functions, we define the matrix-valued differential operators :
 \begin{equation*}
\frac{\partial}{\partial x}  :=  \left(\frac{\partial}{\partial x_{jk}}\right)_{j, k},\qquad 
\frac{\partial}{\partial y}  :=   \left(\frac{\partial}{\partial y_{jk}}\right)_{j, k},\qquad  
\frac{\partial}{\partial z}  :=   \left(\frac{\partial}{\partial x_{jk}} - i \frac{\partial}{\partial y_{jk}}\right)_{j, k},
\end{equation*}We also define: \begin{equation*} \left(\frac{\partial}{\partial z}\right)_{ij}^2 := \sum_k \frac{\partial^2}{\partial z_{ik}\partial z_{kj}},\quad \left(\frac{\partial}{\partial x}\frac{\partial}{\partial y}\right)_{ij} := \sum_k\frac{\partial^2}{\partial x_{ik}\partial y_{kj}}\end{equation*}
\begin{pro}Let functions $f$ satisfying: 
\begin{equation*}\frac{\partial f}{\partial x_{ij}} = \frac{\partial f}{\partial x_{ji}},\quad \frac{\partial f}{\partial y_{ij}} = - \frac{\partial f}{\partial y_{ji}}\qquad \textrm{for all}\quad i, j.\end{equation*}
Then, the infinitesimal generator of a Laguerre process $L(n, m, x)$ is given by: \begin{equation}\mathscr{L} = 2n \,\tr(\Re\left(\frac{\partial}{\partial z}\right)) + 2[\tr(x \Re\left(\frac{\partial}{\partial z}\right)^2) + \tr(y \Im\left(\frac{\partial}{\partial z}\right)^2)]\end{equation}\end{pro}where $\displaystyle \frac{\partial}{\partial z}$ is the operator defined above.

\begin{nota}Using the fact that $x^T = x$, $y^T = -y$ and $\tr(AB) = \tr(BA) = \tr(B^TA^T)$ for any two matrices $A$ and $B$, we can see that 
\begin{equation*}
tr\left(y\frac{\partial}{\partial y}\frac{\partial}{\partial x}\right) = tr\left(\frac{\partial}{\partial x}\frac{\partial}{\partial y}y\right) 
= tr\left(y\frac{\partial}{\partial x}\frac{\partial}{\partial y}\right) \Rightarrow
\left(y \Im(\frac{\partial}{\partial z})^2\right) = 2tr\left(y\frac{\partial}{\partial x}\frac{\partial}{\partial y}\right)\end{equation*}
\end{nota}
\section{Eigenvalues of Laguerre Process}
In this section, we will suppose that $n \geq m$.  The following result is due to K\"onig and O'Connell  (2001), Katori and Tanemura (2004) and Bru in the real case (1989a) : 
\begin{teo}\label{exi}
Let $\lambda_1(t), \cdots, \lambda_m(t)$ denote the eigenvalues of $X_t$. Suppose that at time $t=0$, all the eigenvalues are distinct. Then, the eigenvalue process  $(\lambda_1(t), \hdots, \lambda_m(t))$ satisfies the following stochastic differential system: 
\begin{equation*}
d\lambda_i(t) = 2\sqrt{\lambda_i(t)}\, d\beta_i(t) + 2\left[ n + \sum_{k \neq i}\frac{\lambda_i(t) + \lambda_k(t)}{\lambda_i(t) - \lambda_k(t)}\right]dt \qquad 1 \leq i \leq m,\qquad t < \tau, \end{equation*}where the $(\beta_i)_{1 \leq i \leq m}$ are independent Brownian motions and $\tau$ is defined by \\$\tau := \inf\{t, \lambda_i(t) = \lambda_j(t) \, \textrm{for some}\, (i,j)\}$.   
\end{teo}
  
\begin{nota}With the help of the SDE satisfied by the eigenvalues, we can compute the ones satisfied by both processes $(\tr(X))$ and $(\det(X))$: the former is done. For the second, we find that for $t < T_0 := \inf\{t, \det(X_t) = 0\}$ and for $r \in \re$: 
\begin{eqnarray*}
d(\det(X_t)) & = &  2\det(X_t)\sqrt{\tr(X_t^{-1})}d\nu_t + 2(n-m+1)\det(X_t)\tr(X_t^{-1})dt \\
d(\log(\det(X_t)) & = &  2\sqrt{\tr(X_t^{-1})}d\nu_t + 2(n-m)\tr(X_t^{-1})dt, \\
d(\det(X_t)^r) & = & 2r(\det(X_t))^r \sqrt{\tr(X_t^{-1})} d\nu_t + 2r(n-m+r) (\det(X_t))^r \tr(X_t^{-1})dt  
\end{eqnarray*} so, we can see that for $n=m$, $\log(\det(X)$ is a local martingale and so is $(\det(X))^{m-n}$.
\end{nota}\begin{lem} Take $X_0 \in \df$. Then for $n \geq m$, $X_t \in \df$.\end{lem}
{\it Proof} : In fact, this result is a direct consequence of the fact that for $n=m$, $\log \det(X)$ is a local martingale, and so is $(\det (X) )^{m-n}$. Hence, for $n \geq m$, these two continuous processes tend to infinity when $t \rightarrow T_0$  which is possible only if $T_0 = \infty$, because every continuous local martingale is a time-changed Brownian motion.$\hfill \square$
 \begin{cor}If $\lambda_1(0) > \cdots > \lambda_m(0)$, then, the process $U$ defined by \begin{equation*}U(t) = \frac{1}{\prod_{i < j}(\lambda_i(t) - \lambda_j(t))},  \qquad t < \tau \end{equation*} is a  local martingale.\end{cor}
{\it Proof} : We could follow the proof given by Bru (1989a) or make straightforward computations using the derivatives of the Vandermonde function.
But we prefer use a result from K\"onig and O'Connell (2001): for $n \geq m$,  the eigenvalue process is the $V$-transform (in the Doob sense) of the process obtained from $m$ independent $BESQ(2(n - m +1))$. Thus, if G and $\hat{G}$ denote respectively the infinitesimal generators of these two processes, then, $G(h) = 0$ and, for all $C^2$ function $f$,
\begin{equation*} \hat{G}(f) = \frac{1}{V}G(Vf) \Rightarrow  \hat{G}(U) = \frac{1}{V}G({\bf 1}) =0.  \end{equation*} 
\begin{cor}If at time $t=0$, the eigenvalues of $X$ are distinct, then, they will never collide, i.e. $\tau = \infty$ almost surely. \end{cor}
{\it Proof} : This result follows from the fact that the continuous process $U$ tends to infinity when $t \rightarrow \tau$ which is possible only if $\tau = \infty$ almost surely (We use the same argument as before).$\hfill \square$  

\subsection{Additivity Property}
The proof of this result is similar to the one derived by Bru in the real case (1989): 
\begin{pro}If $(X_t)_{t \geq 0}$ and $(Y_t)_{t \geq 0}$ are two independent Laguerre processes $L(n, m, X_0)$ and $L(p, m, Y_0)$ respectively, then the process $(X_t + Y_t)_{t \geq 0}$ is a Laguerre process $L(n+p, m, X_0 + Y_0)$.\end{pro}

Now, we introduce the Laguerre processes of noninteger dimensions $\delta$.
\section{Laguerre Processes with noninteger dimensions}
Let $X$ be a Laguerre process $L(n, m ,X_0)$ with $n \geq m$. If $X_0 \in \tilde{H_m^+}$, and if $\sqrt{X_t}$ stands for the symmetric matrix square root of $X_t$, it is easy to show that the matrix $O$ defined by $O_t := \sqrt{X_t}^{-1} B_t^{\star}$, where $X_t = B_t^{\star}B_t$, satisfies $O^{\star}O = OO^{\star} = {\it I}_m$. Thus, \begin{equation*}d\gamma_t = O_t dB_t = \sqrt{X_t}^{-1} B_t^{\star}dB_t \end{equation*} is a $m \times m$ complex Brownian matrix. Replacing this expression in (\ref{Henrike}), one obtains :
\begin{equation*} dX_t = \sqrt{X_t}d\gamma_t + d\gamma_t^{\star} \sqrt{X_t} + 2n{\it I_m} dt \end{equation*}

\begin{teo}
If $(B_t)$ is a $m \times m$ complex Brownian matrix, then for every $X_0 \in \df$ and for all $\delta \geq m$, the SDE 
\begin{equation}\label{TIG}
dX_t = \sqrt{X_t}dB_t + dB_t^{\star} \sqrt{X_t} + 2\delta{\it I_m} dt \end{equation} has a unique strong solution in $\df$. Furthermore, if  the eigenvalues are distinct at time $t=0$, then they satisfy the stochastic differential system:
\begin{equation*}d\lambda_i(t) = 2\sqrt{\lambda_i(t)}\, d\beta_i(t) + 2\left[ \delta + \sum_{k \neq i}\frac{\lambda_i(t) + \lambda_k(t)}{\lambda_i(t) - \lambda_k(t)}\right]dt \qquad 1 \leq i \leq m, \end{equation*}where the $(\beta_i)_{1 \leq i \leq m}$ are independent Brownian motions.\end{teo} 
{\it Proof} : The proof of the second part of the theorem is the same as before with $\delta$ instead of $n$. So, we have to prove the first part. Note first that $(\det(X_t)), (\log \det(X_t))$ \\and $(\det(X_t)^r)$ verify the same SDE with $\delta$ instead of $n$. Hence, arguing as before, we can see that $T_0 = \infty$ almost surely. On the other hand, the map $a \mapsto a^{1/2}$ is analytic in $\tilde{H_m^+}$ (see Roger and Williams 1987, p. 134), so, the SDE has a unique strong solution for all $t \geq 0. \hfill \square$  
\begin{defi}Such a process is called the {\it Laguerre process} of dimension $\delta$, size $m$ and initial state $X_0$. It will be denoted by $L(\delta, m, X_0)$.\end{defi}
\begin{nota}
Any process $(X_t)_{t \geq 0}$ solution of (\ref{TIG}) is a diffusion whose infinitesimal \\generator is given by:
\begin{equation*}\mathscr{L} = 2\delta \, \tr(\Re(\frac{\partial}{\partial z})) + 2[\tr(x \Re(\frac{\partial}{\partial z})^2) + \tr(y \Im(\frac{\partial}{\partial z})^2)]\end{equation*}
\end{nota}
\begin{nota}
A simple computation shows that \begin{equation*} d\langle X_{ij}, X_{kl}\rangle_t = 2(X_{il}(t) \delta_{kj} + X_{kj}(t) \delta_{il})dt, \quad
\textrm{for all}\, \,i,j,k,l \in \{1,\hdots ,m\} \end{equation*}
\end{nota}

Now, we focus on both existence and uniqueness when $\delta > m-1$ and $X_0 \in H_m^+$ (see Bru (1989) for the real case).
\subsection{The Process $X^+$}
If $X$ is a Hermitian matrix, let $X^+$ be the Hermitian matrix $\max(X , 0)$. If we denote by $(\lambda_i)$ the eigenvalues of $X$, then  $(\lambda_i^+ = \max(\lambda_i, 0))$ are those of $X^+$. 
\begin{teo}
For all $\delta \in {\mathbb{R}}_+$ and $X_0 = x \in H_m$ , the stochastic differential equation 
\begin{equation} \label{Baya}
dX_t = \sqrt{X_t^+}dB_t + dB_t^{\star}\sqrt{X_t^+} + 2\delta {\it I}_m dt \end{equation} has a solution in $H_m$.\end{teo}
{\it Proof} : The mapping $a \mapsto \sqrt{a^+}$ is continuous on $H_m$. Hence, $X$ exists up to its explosion time (Ikeda and Watanabe 1989, theorem. 2. 3). Furthermore, from 
\begin{equation*}||\sqrt{X^+}||^2 + || \delta {\it I}||^2 \leq {\delta}^2 + ||X||^2 \leq C(1 + ||X||^2),\end{equation*}we can deduce that this explosion time is infinite almost surely (Ikeda and Watanabe 1989, theorem 2. 4).$\hfill \blacksquare$ 
\begin{pro}
If $ \lambda_1(0) > \hdots > \lambda_m(0) \geq 0$, then, for all $t < S := \inf\{t, \lambda_i = \lambda_j \,\textrm{for some (i, j)}\}$, the eigenvalues of $X^+$ verify the following differential system:
\begin{equation*}d\lambda_i(t) = 2\sqrt{\lambda_i^+(t)} d\nu_i(t) + 2\left( \delta + \sum_{k \neq i}\frac{\lambda_i^+(t) + \lambda_k^+(t)}{\lambda_i(t) - \lambda_k(t)}\right)dt, \qquad 1 \leq i \leq m,\end{equation*}\end{pro} 
{\it Proof} : This differential system can be shown in the same way as in theorem \ref{exi} using :
\begin{equation*} \langle dX_{ij},dX_{kl}\rangle_t = 2(X_{il}^+(t) \delta_{kj} + X_{kj}^+(t) \delta_{il})dt, \quad \textrm{for all}\, i,j,k,l \in \{1,\hdots ,m\}  \end{equation*}
\begin{pro}\label{pos}
If $\lambda_1(0) > \hdots > \lambda_m(0) \geq 0$, then, for all $\delta > m-1,\, t > 0, \, \lambda_m(t) \geq 0$ .\end{pro}
{\it Proof} : First, we note that $S = \infty$ almost surely. Indeed,  one can easily show that the process $U$ defined by :
\begin{equation*}
U(\lambda_1(t), \hdots, \lambda_m(t)) = \frac{1}{\prod_{i < j}(\lambda_i(t)- \lambda_j(t))}\end{equation*} is a local martingale. 
For the proof, we follow in the same way as Bru (1989).
 \begin{teo}If $\lambda_1(0) > \hdots > \lambda_m(0) \geq 0$, then, for all $\delta > m-1$, $(\ref{TIG})$ has a unique solution in $H_m^+$ in the sense of probability law. \end{teo}
{\it Proof} : By Proposition $\ref{pos}$, the solution of the SDE (\ref{Baya}) remains positive for all $t > 0$, thus, it is a solution of $(\ref{TIG})$.$\hfill \square$ 
\begin{teo}Let $H_m^+$ be the space of positive Hermitian matrices. Then, whenever the SDE (\ref{TIG}) has a solution in $H_m^+$, for fixed t, its distribution is given by its Laplace transform:\begin{equation}
\label{Lap}
\mathbb{E}_{X_0}(\exp-({\tr uX_t})) = (\det( \it I_m + 2tu))^{-\delta} \exp(-{\tr(X_0(\it I_m + 2tu)^{-1}u)}), \end{equation} for all $u$ in $H_m^+$.\end{teo}
{\it Proof} : For $s \in H_m^+$, let $g(t,s) = \Delta_t^{-\delta} \exp(-V(t,s))$ , where \begin{equation*}\Delta_t = \det(\it I_m + 2ut),\quad W_t =(\it I_m +2ut)^{-1}u,\quad V(t,s) = \tr(sW_t),\end{equation*}First, note that $W \in H_m$. To proceed, we need a lemma :
\begin{lem}
The function $g$ satisfies the heat equation:  $\displaystyle \frac{\partial g}{\partial t} = \mathscr{L} g$ where $\mathscr{L}$ is the infinitesimal generator of $X$.\end{lem}
{\it Proof of the lemma} : If we write $s= x + iy$, then, using the fact that $x$ is symmetric, $y$ is skew-symmetric and $W$ is Hermitian, we can see that $\tr(sW_t) = \tr(xM + iyN)$ where \begin{equation*}M = \frac{W + \overline{W}}{2} \qquad N =  \frac{W - \overline{W}}{2}.\end{equation*}
Observing that $M^T = M$ and $N^T = -N$, we can deduce that $g$ satisfies the conditions of Proposition 1. Besides, \begin{eqnarray*}
\frac{\partial g}{\partial t}  &=&  -g( 2\delta\, \tr(W_t) - 2 \tr (sW_t^2))\\
\tr(y\left(\frac{\partial^2 g}{\partial x\partial y} + \frac{\partial^2 g}{\partial y\partial x}\right)) & = & -ig\, \tr(yW^2),\\
\tr(x\left(\frac{\partial^2 g}{\partial x^2} - \frac{\partial^2 g}{\partial y^2}\right)) & = & g\, \tr(x(M^2 + N^2)) = g\,\tr(xW^2).
\end{eqnarray*}Finally, noting that $\tr(M) = \tr(W)$, we obtain the equality. Now, we consider the process $(Z(t, X_t))$ defined by $Z(t, X_t) = g(t_1 -t, X_t)$ for all $t \leq t_1$ for fixed $t_1$. From the lemma, we deduce that Z is a bounded local martingale and thus is a martingale. So, the result follows from a simple application of the optional stopping theorem.$\hfill \blacksquare$
\begin{cor}\label{samia}
Let $(X_t)_{t \geq 0}$ be a Laguerre process $L(\delta, m, x)$ where $x \in \tilde{H_m}^+$. For $\delta > m-1$, its semi-group is given by the following density: 
\begin{equation*}
p_t^{\delta}(x,y) = \frac{1}{(2t)^{m\delta}\Gamma_m(\delta)}\exp-(\frac{1}{2t}\tr(x+y))\,(\det\,y)^{\delta - m} \f(\delta; \frac{xy}{4t^2}){\bf 1}_{\{y > 0\}}\end{equation*}
with respect to Lebesgue measure $dy = \prod_{p\leq q} dy_{pq}^1\prod_{p<q}dy_{pq}^2$ where $y = y^1+iy^2$ and ${}_0F_1$ is a hypergeometric function of Hermitian matrix argument
(Chikuze 1976, Gross and Richards 1991). \end{cor}
{\it Proof} : In fact, this result can be easily deduced from the case where $\delta = n$ is integer, since, in this case, $X_t$ is a non-central complex Wishart variable $W(n, 2t{\it I}_m, x)$ 
(James 1964) with density given by: 
\begin{equation*}f_t(x,y) =  \frac{1}{(2t)^{mn}\Gamma_m(n)}\exp-(\frac{1}{2t}\tr(x+y))\,(\det\,y)^{n - m} \f(n; \frac{xy}{4t^2}){\bf 1}_{\{y > 0\}}\end{equation*}with respect to $dy$. Hence, taking $\delta$ instead of $n$ and denoting by $W_t$ this new variable (starting from $x$), we can see that : (we will use $|y|$ to denote $\det(y)$)
\begin{align*}
E_x(e^{-\tr uW_t}) & = \frac{1}{(2t)^{m\delta}\Gamma_m(\delta)}e^{-\frac{tr x}{2t}}\int_{y > 0}\exp \left(-\frac{1}{2t}\tr(({\it I} + 2ut)y)\right) |y|^{\delta - m} \f(\delta; \frac{xy}{4t^2})dy
\\ & = \frac{2t^{m\delta}|x|^{-\delta}}{\Gamma_m(\delta)}e^{-\frac{tr x}{2t}}\int_{z > 0}\exp(-2t \,\tr(x^{-\frac{1}{2}}({\it I} + 2ut)x^{-\frac{1}{2}}z) |z|^{\delta - m} \f(\delta; z)dz
\\ & = \exp(-\frac{tr x}{2t}) |{\it I} + 2ut|^{-\delta} \exp\left(\tr(\frac{x}{2t}({\it I} + 2ut)^{-1})\right)
\\ & =   |{\it I} + 2ut|^{-\delta} \exp\left(-\frac{1}{2t}\tr(x({\it I} + 2ut)^{-1}({\it I} + 2ut - {\it I}))\right)
\\ & =   |{\it I} + 2ut|^{-\delta} \exp\left(-\tr (x({\it I} + 2ut)^{-1}u)\right)
\end{align*}which is equal to (\ref{Lap}).$\hfill \blacksquare$
\begin{nota}
In the last proof, we used the change of variables $z = x^{1/2}yx^{1/2}$ which gives $dz = |x|^m dy$. For the second integral, see Faraut and Kor\`anyi (1994), proposition  XV.1.3.
\end{nota}
\begin{nota}
The expression of the semi-group extends continuously to the degenerate case, namely:
\begin{equation*}
p_t^{\delta}(0_m,y) = \frac{1}{(2t)^{m\delta}\Gamma_m(\delta)}\exp-(\frac{\tr(y)}{2t})\,(\det\,y)^{\delta - m} {\bf 1}_{\{y > 0\}}
\end{equation*}where $0_m$ denotes the null matrix.
\end{nota}
\begin{cor}For $\delta > m-1$, the semi-group of eigenvalue process is given by: \begin{equation*}
q_t(x, y)= \frac{V(y)}{V(x)}\det\left(\frac{1}{2t}\left(\frac{y_j}{x_i}\right)^{\nu/2}e^{-\frac{(x_i+y_j)}{2t}} {\it I}_{\nu}(\frac{\sqrt{x_iy_j}}{t})\right)\end{equation*}
where $x = (x_1, \hdots, x_m), y = (y_1,\hdots, y_m)$ so that $x_1 > \hdots > x_m > 0,\, y_1> \hdots > y_m > 0$, $\delta = m + \nu$ such that $\nu > -1$ and ${\it I}_{\nu}$  denotes the modified Bessel function (Lebedev 1972).\end{cor}
{\it Proof} : The expression of the semi-group can be computed using Karlin and MacGregor formula (1959) since, for $\delta > m-1$, the eigenvalue process is the $h$-transform of the process consisted of $m$ independent $BESQ(2(\delta - m +1))$ conditioned never to collide, as stated by K\"onig and O'Connell (2001).
Another proof is given by P\'ech\'e (2003, p. 68). Here, we will deduce the expression of $q_t(x,y)$ from $p_t(x,y)$ following Muirhead (1982), namely, by projection on the unitary group : we will use Weyl integration formula, then give a determinantal representation of hypergeometric functions of two matrix arguments. 
First, we state Weyl integration formula (Faraut 2006) in the complex case: for any Borel function $f$,\begin{equation*} 
\int_{H_m} f(A) dA = C_m\int_{U(m)}\int_{\re^m}f(uau^*)\alpha(du)(V(a))^2da_1\hdots da_m,\end{equation*} where $C_m = \displaystyle \frac{\pi^{m(m-1)}}{\g(m)}$, $U(m)$ is the unitary group, $\alpha$ is the normalized Haar measure on $U(m)$, $a = diag(a_i)$ and $A=uau^*$. Hence, the semi-group of the eigenvalue process is given by (James 1964):
\begin{align*}
q_t(x, y) &= C_m (V(y)^2)\int_{U(m)}p_t(\tilde{x}, u\tilde{y}u^*)\alpha(du) 
\\&  = \frac{C_m (V(y)^2)}{(2t)^{m\delta}\Gamma_m(\delta)}\prod_{i,j=1}^m e^{-\frac{x_i+y_j}{2t}}\,\left(\prod_{i=1}^m y_j\right)^{\delta - m} \int_{U(m)}\f(\delta; \frac{\tilde{x}u\tilde{y}u*}{4t^2})\alpha(du)
\\& = \frac{\pi^{m(m-1)} (V(y)^2)}{(2t)^{m(m+\nu)}\Gamma_m(m)\g(m+\nu)}\prod_{i,j=1}^m e^{-\frac{x_i+y_j}{2t}}\,\left(\prod_{i=1}^m y_j\right)^{\nu} \f(m + \nu; \frac{\tilde{x}}{4t^2}; \tilde{y}),\end{align*}where $\tilde{y} = diag(y_j)$ , $x$ is a positive definite matrix whose eigenvalues are $x_1, \hdots, x_m$, ${}_0F_1$ in the third line is an hypergeometric function with two matrix arguments (Gross and Richards, 1991) and $\delta = m + \nu, \nu > -1$. Next, we need a lemma.

\begin{lem}Let $B, C \in H_m$ and let $(b_i),\,(c_i)$ be respectively their eigenvalues. Then, \begin{align*}
{}_pF_q((m+\mu_i)_{1\leq i \leq p}, &(m+\phi_j)_{1\leq j\leq q}; B, C) =\pi^{\frac{m(m-1)}{2}(p-q-1)} \g(m)\prod_{i=1}^p\frac{(\Gamma(\mu_i+1))^m}{\g(m+\mu_i)} \\&\prod_{j=1}^q\frac{\g(m+\phi_j)}{(\Gamma(\phi_j+1))^m}\quad \frac{\det\left({}_p\mathscr{F}_q((\mu_i+1)_{1\leq i \leq p}, (1+\phi_j)_{1\leq j\leq q};b_lc_f\right)_{l,f}}{V(B)V(C)}   \end{align*} 
for all $\mu_i,\phi_j > -1, 1\leq i \leq p,\, 1\leq j \leq q$.      \end{lem}
{\it Proof} : Recall that the hypergeometric function of two matrix arguments is given by the following series:\begin{equation*}
{}_pF_q((a_i)_{1\leq i\leq p},(e_j)_{1\leq j\leq q}; B, C) =  \sum_{k=0}^{\infty}\sum_{\tau}\frac{\prod_{i=1}^p(a_i)_{\tau}}{\prod_{j=1}^q(e_j)_{\tau}}\frac{C_{\tau}(B)C_{\tau}(C)}{C_{\tau}(I) k!},\end{equation*}
It is well known that \begin{equation*}C_{\tau}(B) = \frac{k! d_{\tau}}{(m)_{\tau}} s_{\tau}(b_1,\hdots, b_m),\end{equation*}where $s_{\tau}$ is the Schur function and $d_{\tau} = s_{\tau}(1,\dots,1)$ is the representation trace or degree (Gross and Richards 1991, Faraut 2006). Substituting in the series, one gets:
\begin{equation*}{}_pF_q((m+\mu_i)_{1\leq i\leq p},(m+\phi_j)_{1\leq j\leq q}; B, C)= \sum_{k=0}^{\infty}\sum_{\tau}\frac{\prod_{i=1}^p(m+\mu_i)_{\tau}}{\prod_{j=1}^q(m+\phi_j)_{\tau}}\frac{s_{\tau}(B)s_{\tau}(C)}{(m)_{\tau}},\end{equation*}
Now, we write: 
\begin{align*}(m+\mu_i)_{\tau} &= \prod_{r=1}^m\frac{\Gamma(\mu_i+m+k_r-r+1)}{\Gamma(\mu_i+m-r+1)}
= \prod_{r=1}^m\frac{\Gamma(\mu_i+1+k_r +\delta_r)}{\Gamma(\mu_i+m-r+1)}
\\& = \pi^{m(m-1)/2}\frac{(\Gamma(\mu_i+1))^m}{\g(m+\mu_i)}\prod_{r=1}^m(\mu_i+1)_{k_r+\delta_r} , \qquad \delta_r=m-r
\end{align*}
Doing the same for each $(m+\phi_j)_{\tau}$ and for $(m)_{\tau}$, we can see that:
\begin{align*}{}_pF_q((m+\mu_i)_{1\leq i\leq p},(m+\phi_j)_{1\leq j\leq q}; &B, C) =
\pi^{\beta}\g(m)\prod_{i=1}^p\frac{(\Gamma(\mu_i+1))^m}{\g(m+\mu_i)} \prod_{j=1}^q\frac{(\Gamma(\phi_j+1))^m}{\g(m+\phi_j)}
\\&\sum_{k=0}^{\infty}\sum_{\tau}\prod_{r=1}^m\left(\frac{\prod_{i=1}^p(\mu_i+1)_{k_r+\delta_r}}{\prod_{j=1}^q(\phi_j+1)_{k_r+\delta_r}}\right)\quad
\frac{s_{\tau}(B)s_{\tau}(C)}{\prod_{r=1}^m(1)_{k_r+\delta_r}}\end{align*}where $\beta = \displaystyle \left(\frac{m(m-1)}{2}\right)(p-q-1)$. To get the desired result, we use the ''Hua formula'' (Faraut 2006) :
\begin{lem} Given an entire function f, i.e. $f(z) = \sum_{k=0}^{\infty}e_k z^k$, then \begin{equation*}
\frac{\det(f(b_ic_j))_{i,j}}{V(B)V(C)} = \sum_{k=0}^{\infty}\sum_{\tau}\left(\prod_{r=1}^m e_{k_r+\delta_r}\right)\,\frac{s_{\tau}(B)s_{\tau}(C)}{s_{\tau}({\it I}_m)}.\end{equation*}\end{lem}  Thus, we get: \begin{align*}
{}_pF_q((m+\mu_i)_{1\leq i \leq p}, &(m+\phi_j)_{1\leq j\leq q}; B, C) =\pi^{\frac{m(m-1)}{2}(p-q-1)}\g(m) \prod_{i=1}^p\frac{\Gamma(\mu_i+1)}{\g(m+\mu_i)} \\&\prod_{j=1}^q\frac{\g(m+\phi_j)}{\Gamma(\phi_j+1)}\quad \frac{\det\left(\sum_{k=0}^{\infty}\frac{\prod_{i=1}^p(\mu_i+1)_k}{\prod_{j=1}^q(\phi_j+1)_k}\, \frac{(b_lc_p)^k}{k!}\right)_{l,p}}{V(B)V(C)} 
\\&=\pi^{\frac{m(m-1)}{2}(p-q-1)}\g(m) \prod_{i=1}^p\frac{\Gamma(\mu_i+1)}{\g(m+\mu_i)} \prod_{j=1}^q\frac{\g(m+\phi_j)}{\Gamma(\phi_j+1)}\\& \frac{\det\left({}_p\mathscr{F}_q((\mu_i+1)_{1\leq i \leq p}, (1+\phi_j)_{1\leq j\leq q}; b_lc_f\right)_{l,f}}{V(B)V(C)} \qquad  \square \end{align*}

For $p=0$ and $q\geq 1$, it reads : \begin{align*}
{}_0F_q((m+\phi_j)_{1\leq j\leq q}; B, C) =\pi^{-\frac{m(m-1)}{2}(q+1)} \g(m)& \prod_{j=1}^q\frac{\g(m+\phi_j)}{(\Gamma(\phi_j+1))^m}\\& \frac{\det\left({}_0\mathscr{F}_q( (1+\phi_j)_{1\leq j\leq q};b_lc_f\right)_{l,f}}{V(B)V(C)},   \end{align*}and similarly,\begin{equation*}
{}_0F_0(B, C) = \frac{\g(m)}{\pi^{\frac{m(m-1)}{2}}}\, \frac{\det(e^{b_lc_f})_{l, f}}{V(B)V(C)}\end{equation*} which can be viewed as Harish-Chandra formula for the "Itzykson-Zuber" integral (Collins 2003). We now proceed to the end of the proof. Taking $p=0,\,q=1,\,B=\displaystyle \frac{\tilde{x}}{4t^2},\,C=\tilde{y}$, we get:\begin{equation*} \f(m + \nu; \frac{\tilde{x}}{4t^2}; \tilde{y}) = \frac{(4t^2)^{m(m-1)/2}\g(m+\nu)\g(m)}{\pi^{m(m-1)}(\Gamma(\nu+1))^m}\,\frac{\det\left({}_0\mathscr{F}_1((\nu+1); \frac{x_iy_j}{4t^2})\right)}{V(x)V(y)}\end{equation*}
The expression of $q_t(x, y)$ follows from a simple computation and from the fact that: \begin{equation*}\frac{{}_0\mathscr{F}_1((\nu+1); x_iy_j/4t^2))_{i,j}}{\Gamma(\nu+1)} = \left(\frac{2t}{\sqrt{x_iy_j}}\right)^{\nu} {\it I}_{\nu}(\frac{\sqrt{x_iy_j}}{t}) \hspace{2.5cm} \blacksquare\end{equation*}
\begin{pro}The measure defined by $\rho(dx) = (\det(x))^{\delta-m} \,dx$ on $\tilde{H_m}^+$ is invariant under the semi-group, i. e, $\rho P_t = \rho$.\end{pro}
{\it Proof}: Denote by $P_t$ the semi-group of Laguerre process $L(\delta, m, x)$ for $\delta > m-1$. Then, we have to show that \begin{equation*}\int_{x>0}P_tf(x)\rho(dx) = \int_{y>0}f(y) \rho(dy)\qquad f \in C_0(\tilde{H_m}^+).\end{equation*} 
This follows by a similar computation and the same arguments as in the proof of corollary \ref{samia}. $\hfill \blacksquare$
\begin{nota}For Wishart processes, it is easy to see that $\mu(dx) := (\det(x))^{\frac{\delta}{2} - \frac{m+1}{2}} {\bf 1}_{\{x>0\}}dx$ is invariant under the semi-group.\end{nota} 

\section{Girsanov Formula and  Absolute-continuity Relations}
The index $\nu > -1$ of a $L(\delta,m,x)$ is defined by $\nu = \delta - m$. In this section, we will discuss in the same way as in \cite{Don} to derive absolute-continuity relations between different indices. 
\subsection{Positive Indices}
Take a matrix-valued Hermitian predictable process H. Let $Q_x^{\delta}$ be the probability law of $L(\delta, m, x)$ for $\delta > m-1$ and $x \in \tilde{H_m^+}$. 
Define: \begin{eqnarray*}
L_t  & = &  \int_0^t \frac{\tr(H_s dB_s + \overline{H_s}\,\overline{dB_s})}{2},\\
\Phi_t & = & \exp{(L_t - \frac{1}{2}\int_0^t \tr(H_s^2)ds)},\end{eqnarray*}where $B$ is a complex Brownian matrix under $Q_x^{\delta}$. We can easily see that  the process $\beta$ defined by $\beta_t = B_t - \int _0^t H_sds$ is a Brownian matrix under the probability 
\begin{equation*}\mathbb{P}_x^{H} |_{\mathscr{F}_t}:= \Phi_t \cdot Q_x^{\delta} |_{\mathscr{F}_t},\end{equation*}
Furthermore, $(X_t)_{t \geq 0}$ is a solution of \begin{equation}\label{Ahmed}
dX_t = \sqrt{X_t}d\beta_t + d\beta^{\star}_t \sqrt{X_t} + (\sqrt{X_t} H_t + H_t \sqrt{X_t} + 2\delta{\it I_m})dt .\end{equation}
For $H_t =  \nu \sqrt{X_t}^{-1}$, (\ref{Ahmed}) becomes 
\begin{equation*}
dX_t = \sqrt{X_t}d\beta_t + d\beta^{\star}_t \sqrt{X_t} + 2(\delta + \nu){\it I_m}dt ,\end{equation*} so that $(X_t)_{t \geq 0}$ is a $L(\delta + \nu, m, x)$ under $\mathbb{P}_x^{H}$. 
Thus, we proved that :
\begin{teo}For $\delta > m-1$,\begin{equation}Q_x^{\delta + \nu}|_{\mathscr{F}_t} = \exp{\left(\frac{\nu}{2}\int_0^t \tr (\sqrt{X_s}^{-1}dB_s + \overline{\sqrt{X_s}^{-1}dB_s})-\frac{\nu^2}{2}\int_0^t \tr({X_s}^{-1})ds\right)}\cdot Q_x^{\delta}|_{\mathscr{F}_t}.\end{equation}\end{teo}
\begin{pro} 
\begin{equation}\label{Monia}
Q_x^{m + \nu}|_{\mathscr{F}_t} = \left(\frac{det (X_t)}{\det(x)}\right)^{\nu/2}\exp{\left(-\frac{\nu^2}{2}\int_0^t \tr({X_s}^{-1})ds\right)}\cdot Q_x^{m}|_{\mathscr{F}_t}.
\end{equation}
\end{pro} 
{\it Proof} : We know that $\nabla_u(\det(u)) = \det(u) u^{-1}$, hence, $\nabla_u(\log(\det(u))) = u^{-1}$. Then, using the fact that for $\delta = m$, $(\log(\det(X_t)))$ is a local martingale, we get from It\^o formula that:
\begin{align*} \log(\det(X_t)) &= \log(\det(X_0)) + \int_0^t \tr({X_s}^{-1}(\sqrt{X_s}dB_s + dB_s^{\star}\sqrt{X_s}))
\\& = \log(\det(X_0)) + \int_0^t \tr({\sqrt{X_s}}^{-1}dB_s + \overline{{\sqrt{X_s}}^{-1}dB_s}).\qquad \qquad \qquad \square \end{align*}
From (\ref{Monia}), it follows that: 
\begin{cor}\begin{align*}Q_x^m(\exp\left(-\frac{\nu^2}{2}\int_0^t \tr(X_s^{-1})ds\right) | X_t = y) &= {\frac{\det(y)}{\det(x)}}^{-\nu/2} \frac{p_t^{m + \nu}(x, y)}{p_t^m(x, y)}
\\& = \frac{\g(m)}{\g(m + \nu)}(\det(z))^{\nu/2}\frac{\f(m + \nu, z)}{\f(m, z)} \\& := \frac{\tilde{\it  I_{\nu}}(z)}{\tilde{\it I_0}(z)},\end{align*}where  $\displaystyle z=\frac{xy}{4t^2}$. \end{cor}
Now, we state the following asymptotic result:
\begin{cor}Let $X$ be a Laguerre process $L(m, m, x)$, then, as $t \rightarrow \infty$:\begin{equation*}\frac{4}{(m\log t)^2} \int _0^t \tr(X_s)^{-1}ds \qquad \overset{\mathscr{L}}{\rightarrow} \qquad {\it T_{1}}(\beta)\end{equation*}where ${\it T_1}$ is the first hitting time of $1$ by a standard Brownian motion $\beta$.
\end{cor}
{\it Proof} : From (\ref{Monia}), we deduce that:\begin{align*}Q_x^m(\exp\left(-\frac{2\nu^2}{(m\log t)^2}\int_0^t \tr(X_s^{-1})ds\right) | X_t = ty) &=  \frac{\g(m)}{\g(m + 2\nu/m\log t)} (\det(xy/4t))^{\nu/m\log t} \\& \frac{\f(m + 2\nu/m\log t, xy/4t^2)}{\f(m, xy/4t^2)}.\end{align*}Noting that $(t^{m})^{-\nu/m\log t} = e^{-\nu}$, and since both hypergeometric functions converge to $1$ as $t \rightarrow \infty$, we obtain:
\begin{equation*}Q_x^m\left(\exp(-\frac{2\nu^2}{(m\log t)^2}\int_0^t \tr(X_s^{-1})ds | X_t = ty\right) \qquad \overset{t \rightarrow \infty}{\longrightarrow} \qquad e^{-\nu}
\end{equation*}
Then, since\begin{align*}\lim_{t \rightarrow \infty}t^{m^2}p_t^m(x, 2y)&= \lim_{t \rightarrow \infty}\frac{e^{-\tr(x)/2t}}{\g(m)}e^{-\tr(y)} \f(m, \frac{xy}{2t})
\\& = \frac{e^{-\tr(y)}}{\g(m)}\end{align*}we get:
\begin{align*}Q_x^m&(\exp\left(-\frac{2\nu^2}{(m\log t)^2}\int_0^t \tr(X_s^{-1})ds\right)) \\&= \int_{y > 0} Q_x^m\left(\exp(-\frac{2\nu^2}{(m\log t)^2}\int_0^t \tr(X_s^{-1})ds | X_t = y\right)p_t^m(x, y) dy \\&= \int_{y > 0} Q_x^m\left(\exp(-\frac{2\nu^2}{(m\log t)^2}\int_0^t \tr(X_s^{-1})ds | X_t = ty\right) t^{m^2}p_t^m(x, ty) dy\\& \overset{t \rightarrow \infty}{\longrightarrow} 
e^{-\nu},
\end{align*} by dominated convergence Theorem.$\hfill \blacksquare$

\subsection{Negative Indices}
Take $0 < a \leq \det(x)$. The same computation as in parag. $5. 1$ with $H_t = -\nu \sqrt{X_t}^{-1}, 0 < \nu < 1,$ shows that 
\begin{equation*}
Q_x^{m - \nu}|_{\mathscr{F}_{t \wedge T_a}} = \left(\frac{\det(x)}{\det(X_{t \wedge T_a})}\right)^{\nu/2}\exp{\left(-\frac{\nu^2}{2}\int_0^{t \wedge T_a} \tr({X_s}^{-1})ds\right)}Q_x^{m}|_{\mathscr{F}_{t \wedge T_a}}\end{equation*}where $T_a := \inf\{t, \det(X_t) = a\}$. Letting $a \rightarrow 0$ and using the fact that $T_0 = \infty$ a.s under $Q_x^m$, we get :
\begin{align*}Q_x^{m - \nu}|_{\mathscr{F}_{t \wedge T_0}} &= \left(\frac{\det(x)}{\det(X_t)}\right)^{\nu/2}\exp{\left(
\frac{\nu^2}{2}\int_0^t\tr({X_s}^{-1})ds\right)}Q_x^{m}|_{\mathscr{F}_t}\\& =  \left(\frac{\det(x)}{\det(X_t)}\right)^{\nu}Q_x^{m + \nu}|_{\mathscr{F}_t}
\end{align*}
\begin{pro} For all $t > 0$ and $0 < \nu < 1$,
\begin{equation*}
Q_x^{m - \nu}( T_0 > t) =  \frac{\g(m)}{\g(m+\nu)}\det(\frac{x}{2t})^{\nu}{}_1F_1(\nu,m+\nu,-\frac{x}{2t})
\end{equation*}
\end{pro}
{\it Proof} : From the absolute-continuity relation above, we deduce that :
\begin{equation*}Q_x^{m-\nu} (T_0 > t) = Q_x^{m+\nu}\left(\left(\frac{\det(x)}{\det(X_t)}\right)^{\nu}\right),\end{equation*}On the other hand, using the expression of the semi-group, 
one has :
\begin{align*}Q_x^{\delta}(\det(X_t)^s) &= (2t)^{ms}\frac{\g(s + \delta)}{\g(\delta)}{}_1F_1(-s; \delta;
-\frac{x}{2t}) \\& =(2t)^{ms}\frac{\g(s + \delta)}{\g(\delta)}\exp(-\tr(\frac{x}{2t})){}_1F_1(\delta +s; \delta;\frac{x}{2t}) \end{align*}by Kummer relation (cf Th 7. 4. 3 in Muirhead 1982).
Taking $s= -\nu$, we are done.$\hfill \blacksquare$   

\section{Generalized Hartman-Watson law}
Henceforth, we will write $\mathscr{F}$ to denote one-dimensional hypergeometric functions. We define the {\it generalized Hartman-Watson} law as the law of 
\begin{equation*}
\int_0^t \tr(X_s^{-1})ds \qquad \textrm{under} \qquad Q_x^m( \cdot  | X_t = y).\end{equation*} Its Laplace transform is given by:
\begin{equation}\label{G}
Q_x^m(\exp\left(\frac{-\nu^2}{2}\int_0^t \tr(X_s^{-1})ds\right)|X_t =y) = \frac{\g(m)}{\g(m + \nu)}\det(z)^{\nu/2}\frac{\f(m + \nu,z)}{\f(m,z)}  \end{equation}
 $z=xy/4t^2$. Recall that for $m=1$, this is the well-known Hartman-Watson law and that its density was computed by Yor (1980). Here, we will investigate the case $m=2$. The Gross and Richards formula is written for $p=0$ and $q=1$ :
\begin{equation*}
\f(m + \nu,z)= \frac{\det(z_i^{m-j}\h(m + \nu -j +1,z_i))}{V(z)},
\end{equation*}
where $(z_i)$ denote the eigenvalues of $z$ and $V(z) = \prod_{i < j}(z_i - z_j)$ is the Vandermonde determinant . Noting that $\g(m + \nu) = \prod_{j=1}^m\Gamma(m + \nu -j +1)$,
 then :
\begin{equation*}
(\ref{G}) = \frac{\det(z_i^{(m-j)/2}{\it I}_{m + \nu -j}(2\sqrt{z_i}))}{\det(z_i^{(m-j)/2}{\it I}_{m -j}(2\sqrt{z_i}))}
\end{equation*}Without loss of generality, we will take $t=1$.
\begin{pro} For $m = 2$, let $\lambda_1 > \lambda_2$ be the eigenvalues of $\sqrt{xy}$. Then, the density of the generalized Hartman-Watson law is given by :
\begin{equation*}
f(v) = \frac{\sqrt{\lambda_1\lambda_2}v}{p\pi \sqrt{2\pi v^3}}\frac{\int_0^1\int_0^{\infty}z\sinh(p\sqrt{1 - z^2})e^{-2\sqrt{\lambda_1\lambda_2}z \cosh y} e^{-\frac{2(y^2-\pi^2)}{v}}(\sinh y)\sin(\frac{4\pi y}{v})dzdy}{\int_0^1\int_0^1\frac{u\cosh(pu\sqrt{1-x^2})}{\sqrt{1-x^2}}\it I_0(2\sqrt{\lambda_1\lambda_2}ux)dudx},
\end{equation*}
for $v > 0$, where $p = \lambda_1 - \lambda_2$. Furthermore, if $\lambda_1 = \lambda_2 := \lambda$, then: 
\begin{equation*}f(v) = \frac{4\lambda v e^{\frac{2\pi^2}{v}}}{\pi^2 \sqrt{2\pi v^3}}\frac{\int_0^{\infty}g(y) e^{-\frac{2y^2}{v}}(\sinh y)\sin(\frac{4\pi y}{v})dy}{{}_1\mathscr{F}_2(\frac{1}{2}; 1; 2; \lambda^2)},\end{equation*}where \begin{equation*}g(y) = \frac{1}{3} + \frac{\pi}{2}\frac{{\it I}_2(2\lambda \cosh y) + \mathbf{L}_2(2\lambda \cosh y)}{2\lambda \cosh y},\end{equation*}and $\mathbf{L}_2$ is the Struve function (Gradshteyn and Ryzhik, 1994).\end{pro} 
{\it Proof} : For $m=2$ , (\ref{G}) becomes:
\begin{equation*}
(\ref{G}) = \frac{\lambda_1{\it I}_{\nu + 1}(\lambda_1){\it I}_{\nu}(\lambda_2) - \lambda_2{\it I}_{\nu + 1}(\lambda_2){\it I}_{\nu}(\lambda_1)}{ \lambda_1\it I_1(\lambda_1)\it I_0(\lambda_2) - \lambda_2\it I_1(\lambda_2)\it I_0(\lambda_1)},
\end{equation*} so, using the  integral representations below (Brychkov, Marichev, Prudnikov 1986, p. 46) :
\begin{equation*}
x(a\it I_{\nu+1}(ax)\it I_{\nu}(bx) - b\it I_{\nu+1}(bx)\it I_{\nu}(ax)) = (a^2 - b^2)\int_0^x u\it I_{\nu}(au)\it I_{\nu}(bu)du\end{equation*} with $x=1,a=\lambda_1,b=\lambda_2$, and 
(Gradshteyn and Ryzhik 1994, p. 734):
\begin{equation*}
\frac{\pi}{2}\it I_{\nu}(\frac{a}{2}(\sqrt{b^2+c^2}+b))\it I_{\nu}(\frac{a}{2}(\sqrt{b^2+c^2}-b))=\int_0^a\frac{\cosh(b\sqrt{a^2-x^2})}{\sqrt{a^2-x^2}}\it I_{2\nu}(cx)dx
\end{equation*}
where $a>0, \Re(\nu)>-1$, with $a=1,b=(\lambda_1-\lambda_2)u := pu$ et $c=2\sqrt{\lambda_1\lambda_2}u$, the numerator of (\ref{G}) is then equal to:
\begin{equation*}
\frac{2}{\pi}(\lambda_1^2 - \lambda_2^2)\int_0^1\int_0^1\frac{u\cosh(pu\sqrt{1-x^2})}{\sqrt{1-x^2}}\it I_{2\nu}(2\sqrt{\lambda_1\lambda_2}ux)dudx.\end{equation*}
Taking $\nu = 0$, the denominator is then equal to: 
\begin{equation*}
\frac{2}{\pi}(\lambda_1^2 - \lambda_2^2)\int_0^1\int_0^1\frac{u\cosh(pu\sqrt{1-x^2})}{\sqrt{1-x^2}}\it I_0(2\sqrt{\lambda_1\lambda_2}ux)dudx.\end{equation*}
Thus, (\ref{G}) becomes:
\begin{equation*}\frac{\int_0^1\int_0^1\frac{u\cosh(pu\sqrt{1-x^2})}{\sqrt{1-x^2}}\it I_{2\nu}(2\sqrt{\lambda_1\lambda_2}ux)dudx}{\int_0^1\int_0^1\frac{u\cosh(pu\sqrt{1-x^2})}{\sqrt{1-x^2}}\it I_0(2\sqrt{\lambda_1\lambda_2}ux)dudx}\end{equation*}
Now, we only have to use the integral representation of $\it I_{2\nu}$ (Yor 1980):
\begin{align*}\it I_{2\nu}(2\sqrt{\lambda_1\lambda_2}ux)&=\frac{1}{2i\pi}\int_{C}e^{2\sqrt{\lambda_1\lambda_2}ux \cosh \omega}e^{-2\nu\omega}d\omega
\\&=\frac{1}{2i\pi}\int_{C}e^{2\sqrt{\lambda_1\lambda_2}ux\cosh \omega}\int_0^{\infty}\frac{2\omega e^{-v\nu^2/2}}{(2\pi v^3)^{1/2}}e^{-\frac{2\omega^2}{v}}dvd\omega
\end{align*}where $C$ is the contour indicated in Yor (1980), hence, the density function is given by:
\begin{equation*}f(v) = \frac{1}{i\pi \sqrt{2\pi v^3}}\frac{\int_0^1\int_0^1\int_{C}u\omega\frac{\cosh(pu\sqrt{1-x^2})}{\sqrt{1-x^2}}e^{2\sqrt{\lambda_1\lambda_2}ux \cosh \omega} e^{-\frac{2\omega^2}{v}}dudxd\omega}{\int_0^1\int_0^1\frac{u\cosh(pu\sqrt{1-x^2})}{\sqrt{1-x^2}}\it I_0(2\sqrt{\lambda_1\lambda_2}ux)dudx}{\bf 1}_{\{v > 0\}}
\end{equation*}
We can simplify this expression by integrating over C to see that the numerator is equal to (Yor 1980): 
\begin{equation*}\frac{\sqrt{\lambda_1\lambda_2}v}{\pi \sqrt{2\pi v^3}}
\int_0^1\int_0^1\int_0^{\infty}u^2x \frac{\cosh(pu\sqrt{1-x^2})}{\sqrt{1-x^2}}e^{-2\sqrt{\lambda_1\lambda_2}ux \cosh y} e^{-\frac{-2(y^2-\pi^2)}{v}}(\sinh y)\sin(\frac{4\pi y}{v})dudxdy
\end{equation*}
Setting $z = ux$, The numerator is written
 \begin{equation*}\frac{\sqrt{\lambda_1\lambda_2}v}{\pi \sqrt{2\pi v^3}}
 \int_0^1\int_0^u\int_0^{\infty}z \frac{u\cosh(p\sqrt{u^2-z^2})}{\sqrt{u^2-z^2}}e^{-2\sqrt{\lambda_1\lambda_2}z \cosh y} e^{-\frac{2(y^2-\pi^2)}{v}}(\sinh y)\sin(\frac{4\pi y}{v})dudzdy,\end{equation*}
 that we can integrate with respect to $u$ to get \begin{equation*}\frac{\sqrt{\lambda_1\lambda_2}v}{p\pi \sqrt{2\pi v^3}}\int_0^1\int_0^{\infty}z\sinh(p\sqrt{1 - z^2})e^{-2\sqrt{\lambda_1\lambda_2}z \cosh y} e^{-\frac{2(y^2-\pi^2)}{v}}(\sinh y)\sin(\frac{4\pi y}{v})dzdy.\end{equation*}
 
Now, we prove the second part. In this case, $p = 0$ and we have to evaluate :
\begin{equation*}\frac{\lambda ve^{\frac{2\pi^2}{v}}}{\pi \sqrt{2\pi v^3}}\frac{\int_0^1\int_0^1\int_0^{\infty}\frac{u^2x}{\sqrt{1-x^2}}e^{-2\lambda ux \cosh y}
e^{-\frac{2y^2}{v}}(\sinh y)\sin(\frac{4\pi y}{v})dudxdy}{\int_0^1\int_0^1\frac{u\it I_0(2\lambda ux)}{\sqrt{1-x^2}}dudx}\end{equation*} 
Setting $z = ux$, the numerator reads :
\begin{equation*}\frac{\lambda v e^{\frac{2\pi^2}{v}}}{\pi \sqrt{2\pi v^3}}\int_0^1\int_0^{\infty}z \sqrt{1-z^2}e^{-2\lambda z \cosh y} e^{-\frac{2y^2}{v}}(shy)sin(\frac{4\pi y}{v})dzdy,\end{equation*}
Integration with respect to $z$ yields (Gradshteyn and Ryzhik 1994, p. 369):
\begin{equation*}\frac{\lambda v e^{\frac{2\pi^2}{v}}}{\pi \sqrt{2\pi v^3}}\int_0^{\infty}g(y) e^{-\frac{2y^2}{v}}(shy)sin(\frac{4\pi y}{v})dy\end{equation*}
For the denominator, we use the fact that $\frac{d}{dz}(z\it I_1(z)) = z\it I_0(z)$, which yields:
\begin{equation*}
\int_0^1\int_0^1\frac{u\it I_0(2\lambda ux)}{\sqrt{1-x^2}}dudx=\int_0^1\frac{\it I_1(2\lambda x)}{2\lambda x\sqrt{1-x^2}}dx\end{equation*}
Then, the following formula
\begin{align}\label{Mehdi}
\int_0^ax^{\alpha - 1}(a^2-x^2)^{\beta-1}\it I_{\nu}(cx)dx &= 2^{-\nu -1}a^{2\beta + \alpha + \nu -2}c^{\nu}\frac{\Gamma(\beta)\Gamma((\alpha + \nu) /2)}{\Gamma(\beta+(\alpha+\nu)/2)\Gamma(\nu+1)}
\\&{}_1\mathscr{F}_2(\frac{\alpha+\nu}{2}; \beta +\frac{\alpha+\nu}{2}; \nu+1; \frac{a^2c^2}{4})\nonumber
\end{align}taken with $\alpha =0, a=1,\beta=1/2,c=2\lambda,\nu=1$ gives:\begin{equation*}
\int_0^1\frac{\it I_1(2\lambda x)}{2\lambda x\sqrt{1-x^2}}dx=\frac{\pi}{4}{}_1\mathscr{F}_2(\frac{1}{2}; 1; 2; \lambda^2)\end{equation*}

We can proceed differently :  let $\lambda_1 = \lambda_2 + h$ then (\ref{G}) reads:
\begin{equation*}\frac{((\lambda_2+h){\it I_{\nu+1}}(\lambda_2+h){\it I_{\nu}}(\lambda_2) - \lambda_2{\it I_{\nu+1}}(\lambda_2){\it I_{\nu}}(\lambda_2+h))/h}{((\lambda_2+h){\it I_{1}}(\lambda_2+h){\it I_{0}}(\lambda_2) - \lambda_2{\it I_{1}}(\lambda_2){\it I_{0}}(\lambda_2+h))/h}.\end{equation*}
Next, we let $h \rightarrow 0$. As usual, we first compute the numerator and then take $\nu = 0$. To do this, we shall evaluate :\begin{eqnarray*}
A & = & \lim_{h\rightarrow 0}\frac{(\lambda_2+h){\it I_{\nu+1}}(\lambda_2+h) - \lambda_2 {\it I_{\nu+1}}(\lambda_2)}{h} \\
B & = & \lim_{h\rightarrow 0}\frac{{\it I_{\nu}}(\lambda_2+h) - {\it I_{\nu}}(\lambda_2)}{h}  
\end{eqnarray*} which are equal respectively to $\displaystyle \frac{d}{dx}(x{\it I_{\nu+1}}(x))$ and $\displaystyle \frac{d}{dx}({\it I_{\nu}}(x))$ taken for $x=\lambda = \lambda_1 = \lambda_2$. Using the differentiation formula  $\displaystyle \frac{d}{dx}(x^{\nu} {\it I_{\nu}}(x)) = x^{\nu} {\it I_{\nu-1}}(x)$ (Lebedev 1972, p. 110), we get:
\begin{eqnarray*}
 \frac{d}{dx}(x{\it I_{\nu+1}}(x))  =   -\nu {\it I_{\nu+1}}(x) + x {\it I_{\nu}}(x), \quad
 \frac{d}{dx}({\it I_{\nu}}(x))  =  -\frac{\nu}{x} {\it I_{\nu}}(x) + {\it I_{\nu-1}}(x), 
\end{eqnarray*}  
thus: \begin{eqnarray*}
N&=& {\it I_{\nu}}(\lambda)( -\nu {\it I_{\nu+1}}(\lambda) + \lambda {\it I_{\nu}}(\lambda)) - \lambda {\it I_{\nu+1}}(\lambda)( -\frac{\nu}{\lambda} {\it I_{\nu}}(\lambda) + {\it I_{\nu-1}}(\lambda)) \\
&=& \lambda({\it I_{\nu}}^2(\lambda) - {\it I_{\nu+1}}(\lambda) {\it I_{\nu-1}}(\lambda))\\
A &=&  \frac{{\it I_{\nu}}^2(\lambda) - {\it I_{\nu+1}}(\lambda) {\it I_{\nu-1}}(\lambda)}{{\it I_{0}}^2(\lambda) - {\it I_{1}}(\lambda) {\it I_{-1}}(\lambda)}
\end{eqnarray*} 
Using the integral representation below (Gradshteyn and Ryzhik 1994, p. 757):\begin{equation*}
{\it I_{\mu}}(z) {\it I_{\nu}}(z) = \frac{2}{\pi} \int_0^{\pi/2} \cos((\mu - \nu)\theta) {\it I_{\mu+\nu}}(2z\cos \theta) d\theta ,\qquad \Re(\mu+\nu) > -1.\end{equation*} 
the numerator is written as  :
\begin{align*}
N= \frac{2}{\pi} \int_0^{\pi/2} (1-\cos2\theta){\it I_{2\nu}}(2\lambda \cos \theta) d\theta
&=\frac{4}{\pi} \int_0^{\pi/2} (\sin^2\theta) {\it I_{2\nu}}(2\lambda \cos \theta)d\theta
\\&=\frac{4}{\pi} \int_0^{1} \sqrt{1 - r^2}{\it I_{2\nu}}(2\lambda r) dr,\end{align*}
Thus, using (\ref{Mehdi}), the denominator is equal to  
\begin{equation*}
D= \frac{4}{\pi} \int_0^{1} \sqrt{1 - r^2}{\it I_0}(2\lambda r) dr = \frac{\pi}{4}{}_1\mathscr{F}_2(\frac{1}{2}; 1; 2; \lambda^2)\end{equation*}
Finally,  the integral representation of ${\it I_{\nu}}$ gives :  
\begin{align*}
f(u)&= \frac{\lambda u e^{2\pi^2/u}}{\pi\sqrt{2\pi u^3}}\displaystyle \frac{\int_0^{\infty} e^{-2y^2/u}\sinh(y)\sin\left(\displaystyle \frac{4\pi y}{u}\right)\int_0^1r\sqrt{1-r^2}e^{-2\lambda r\cosh y}dr \,du}{\int_0^1 \sqrt{1 - r^2}{\it I_0}(2\lambda r) dr}
\\& = \frac{\lambda u e^{2\pi^2/u}}{\pi\sqrt{2\pi u^3}}\displaystyle \frac{\int_0^{\infty}g(y) e^{-2y^2/u}\sinh(y)\sin\left(\displaystyle \frac{4\pi y}{u}\right)du}{\int_0^1 \sqrt{1 - r^2}{\it I_0}(2\lambda r) dr},\end{align*}
  
\section{The Law of $T_0$}
Recall that: For $0 < \nu < 1$,\begin{equation*}
Q_x^{m-\nu}(T_0 > t) = \frac{\g(m)}{\g(m+\nu)}\det(\frac{x}{2t})^{\nu}{}_1F_1(\nu,m+\nu,-\frac{x}{2t})\end{equation*}
\begin{pro} Let $m=2$ and $\lambda_1 > \lambda_2$ be the eigenvalues of $x$. The density of $ S_0 := 1/(2T_0)$ under $Q_x^{m-\nu}$ is given by: 
\begin{equation*}
f(u)  = \frac{(\lambda_1\lambda_2)^{\nu}u^{2\nu-2}e^{-(\lambda_1 + \lambda_2)u}}{\Gamma(\nu+1)\Gamma(\nu)}\, \frac{\fc(2,\nu+1,\lambda_1u) - \fc(2,\nu+1,\lambda_2u)}{(\lambda_1 - \lambda_2)}\end{equation*}\end{pro}
\begin{cor}If $\lambda_1 = \lambda_2 := \lambda$, the density is written:\begin{equation*} f(u) =  \frac{2\lambda^{2\nu}u^{2\nu-1}e^{-\lambda u}}{\Gamma(\nu+2)\Gamma(\nu)}\fc(\nu-1,\nu+2,-\lambda u)\end{equation*}\end{cor}
{\it Proof} : Recall first that when $m=1$, $S_0 \overset{\mathscr{L}}{=} \gamma_{\nu}/x$, where $\gamma_{\nu}$ is a Gamma variable with density $r^{\nu-1}e^{-r}dr$. With the help of the Gross-Richards formula, it follows that for $m=2$,\begin{align*}
Q_x^{m-\nu}(S_0 \leq u) &= \frac{(\lambda_1\lambda_2)^{\nu}}{(\lambda_1 - \lambda_2)\Gamma_2(\nu + 2)}u^{2\nu}(\lambda_1\fc(\nu,\nu+2,-\lambda_1u)\fc(\nu-1,\nu+1,-\lambda_2u) \\& - \lambda_2\fc(\nu,\nu+2,-\lambda_2u)\fc(\nu-1,\nu+1,-\lambda_1u)),\end{align*} where $S_0:=1/(2T_0)$.
This is a $C^{\infty}$ function in $u$. Hence, we will compute its derivative to get the density. Recall that :\begin{equation*}\frac{d}{dz} \, \fc(a,b,z) = \frac{a}{b} \, \fc(a+1,b+1,z),\end{equation*}
thus : 
\begin{equation*}f(u) = \frac{d}{du}Q_x^{m-\nu}(S_0 \leq u) = K(\nu,\lambda_1,\lambda_2)u^{2\nu-1}(A - B)\end{equation*}
where \begin{eqnarray*}
K(\nu,\lambda_1,\lambda_2) & = & \frac{(\lambda_1\lambda_2)^{\nu}}{\Gamma_2(\nu+2)(\lambda_1 - \lambda_2)} \\
A & = & 2\nu((\lambda_1\fc(\nu,\nu+2,-\lambda_1u)\fc(\nu-1,\nu+1,-\lambda_2u) \\&-& \lambda_2\fc(\nu,\nu+2,-\lambda_2u) \fc(\nu-1,\nu+1,-\lambda_1u))\\ 
B & = & \frac{\nu}{\nu+2}((\lambda_1^2u\fc(\nu+1,\nu+3,-\lambda_1u)\fc(\nu-1,\nu+1,-\lambda_2u)\\& - & \lambda_2^2u\fc(\nu+1,\nu+3,-\lambda_2u) \fc(\nu-1,\nu+1,-\lambda_1u)).
\end{eqnarray*}Then, we use the contiguous relation :\begin{equation*}
b\, \fc(a,b,z) - b\, \fc(a-1,b,z) = z\,\fc(a,b+1,z)\end{equation*}   to see that \begin{equation*}
\lambda_1u\fc(\nu+1,\nu+3,-\lambda_1u) = (\nu+2)(\fc(\nu,\nu+2,-\lambda_1u)- \fc(\nu+1,\nu+2,-\lambda_1u)) \\
\end{equation*}
\begin{equation*}
\lambda_2u\fc(\nu+1,\nu+3,-\lambda_2u) = (\nu+2)(\fc(\nu,\nu+2,-\lambda_2u)- \fc(\nu+1,\nu+2,-\lambda_2u)) \\
\end{equation*} implies that:
\begin{equation*}f(u) = K_1(\nu,\lambda_1,\lambda_2)u^{2\nu-1}(C+D-E-F)\end{equation*} where
\begin{eqnarray*}
K_1(\nu,\lambda_1,\lambda_2) & = & \frac{\nu(\lambda_1\lambda_2)^{\nu}}{\Gamma_2(\nu+2)(\lambda_1 - \lambda_2)}   \\
C & = & \lambda_1\fc(\nu,\nu+2,-\lambda_1u)\fc(\nu-1,\nu+1,-\lambda_2u)\\
D & = & \lambda_1\fc(\nu+1,\nu+2,-\lambda_1u)\fc(\nu-1,\nu+1,-\lambda_2u) \\
E & = & \lambda_2\fc(\nu,\nu+2,-\lambda_2u)\fc(\nu-1,\nu+1,-\lambda_1u) \\
F & = &  \lambda_2\fc(\nu+1,\nu+2,-\lambda_2u)\fc(\nu-1,\nu+1,-\lambda_1u),
\end{eqnarray*}Applying again the above contiguous relation yields:
\begin{equation*}\lambda_1u\fc(\nu+1,\nu+2,-\lambda_1u) = (\nu+1)(\fc(\nu,\nu+1,-\lambda_1u)- \fc(\nu+1,\nu+1,-\lambda_1u))
\end{equation*}
\begin{equation*}
\lambda_2u\fc(\nu+1,\nu+2,-\lambda_2u) = (\nu+1)(\fc(\nu,\nu+1,-\lambda_2u)- \fc(\nu+1,\nu+1,-\lambda_2u))
\end{equation*}
\begin{equation*}
\lambda_2u\fc(\nu,\nu+2,-\lambda_2u) = (\nu+1)(\fc(\nu-1,\nu+1,-\lambda_2u)- \fc(\nu,\nu+1,-\lambda_2u))
\end{equation*}
\begin{equation*}
\lambda_1u\fc(\nu,\nu+2,-\lambda_1u) = (\nu+1)(\fc(\nu-1,\nu+1,-\lambda_1u)- \fc(\nu,\nu+1,-\lambda_1u))
\end{equation*}
Replacing in the expression of $f$, we obtain \begin{equation*}
f(u) = K_2(\nu,\lambda_1,\lambda_2) u^{2\nu-2}(G-H),\end{equation*}where
\begin{eqnarray*}
K_2(\nu,\lambda_1,\lambda_2) & = & \frac{\nu(\nu+1)(\lambda_1\lambda_2)^{\nu}}{\Gamma_2(\nu+2)(\lambda_1 - \lambda_2)} \\
G & = & \fc(\nu+1,\nu+1,-\lambda_2u)\fc(\nu-1,\nu+1,-\lambda_1u)\\ 
H & = & \fc(\nu+1,\nu+1,-\lambda_1u)\fc(\nu-1,\nu+1,-\lambda_2u)
\end{eqnarray*}Eventually, writing \begin{eqnarray*}
\Gamma_2(\nu+2) & = & \Gamma(\nu+2)\Gamma(\nu+1) =  \nu(\nu+1)\Gamma(\nu+1)\Gamma(\nu)\\
\fc(a,a,z) & = & e^{-z} \qquad    \fc(a,b,-z)  =  e^{-z}\fc(b-a,b,z), 
\end{eqnarray*}we get \begin{equation*}
f(u)  = \frac{(\lambda_1\lambda_2)^{\nu}u^{2\nu-2}e^{-(\lambda_1 + \lambda_2)u}}{\Gamma(\nu+1)\Gamma(\nu)}\, \frac{\fc(2,\nu+1,\lambda_1u) - \fc(2,\nu+1,\lambda_2u)}{(\lambda_1 - \lambda_2)}\end{equation*}
Th case $\lambda_1 = \lambda_2$ is treated in the same way as before (for the Hartman-Watson law). In fact , writing $\lambda_1 = \lambda_2 + h$ and letting $h \rightarrow 0$, we see that the density is given by :\begin{align*}
f(u)  =  \frac{\lambda^{2\nu}u^{2\nu-2}e^{-2\lambda u}}{\Gamma(\nu+1)\Gamma(\nu)}\frac{d}{d\lambda} \, \fc(2,\nu+1,\lambda u)
       & =  \frac{2\lambda^{2\nu}u^{2\nu-1}e^{-2\lambda u}}{\Gamma(\nu+2)\Gamma(\nu)}\fc(3,\nu+2,\lambda u) 
     \\& =  \frac{2\lambda^{2\nu}u^{2\nu-1}e^{-\lambda u}}{\Gamma(\nu+2)\Gamma(\nu)}\fc(\nu-1,\nu+2,-\lambda u) \blacksquare\end{align*}

\section{Conclusion}
The Gross-Richards formula has been the main ingredient in this paper, since it enables us to express more explicitly the special functions of matrix argument. The case $m=3$ can be treated in the same way, but computation becomes too complicated. So, if we want to deal with the general case, it will be convenient to find a more explicit formula. Indeed, Schur functions can be expressed as polynomials in the elementary symmetric functions $e_r$ or as polynomials in the completely symmetric functions $h_r$. More precisely, we have :
\begin{eqnarray*}
s_{\lambda} & = & \det(e_{\lambda_i - i +j}) \qquad  1 \leq i, j \leq n\\
s_{\lambda} & = & \det(h_{\lambda'_i - i + j}) \qquad 1 \leq i, j \leq n  
\end{eqnarray*}where $\lambda$ is a partition of length $\leq n$, and $\lambda'$ is the conjugate of $\lambda$ (Macdonald 1995). So, using these two identities, can we improve our  results?
\section{Appendix: special functions}
\subsection{The hypergeometric series}
The multivariate hypergeometric functions were studied by Muirhead (1982) in the real symmetric case, Chikuze (1976) for the complex Hermitian case and Faraut and Kor\`anyi (1994) in a more general setting. For Hermitian matrix argument, they are defined by: \begin{equation*}
{}_pF_q((a_i)_{1\leq i \leq p}, (b_j)_{1\leq j \leq q}; X) = \sum_{k \geq 0}\sum_{\tau \bot k}\frac{(a_1)_{\tau}\cdots(a_p)_{\tau}}{(b_1)_{\tau}\cdots(b_q)_{\tau}}\frac{C_{\tau}(X)}{k!}
\end{equation*}
where $\tau = (k_1, \hdots, k_m)$ is a partition of weight $k$ and length $m$ such that $k_1 \geq \hdots \geq k_m$, 
$(a)_{\tau}$ is the generalised Pochammer symbol defined by: \begin{align*}
(a)_{\tau} = \prod_{i=1}^m\frac{\Gamma(a+k_i - i + 1)}{\Gamma(a - i + 1)} , \qquad \tau=(k_1, \dots, k_m)\end{align*} 
and $C_{\tau}(X)$ is the zonal polynomial of $X$ such that : \begin{equation*}(\tr(X))^k = \sum_{\tau \bot k}C_{\tau}(X)\end{equation*}
Several normalizations for this polynomial exist in the litterature but we consider this one. This polynomial is symmetric, homogeneous, of degree $k$ in the eigenvalues of $X$ and is an eigenfunction of the following differential operator :
\begin{equation*} \Delta_{X} = \sum_{i=1}^mx_i^2 \frac{\partial^2}{\partial x_i^2} + 2\sum_{i=1}^m\sum_{1\leq k \neq i \leq m}\frac{x_i^2}{x_i - x_k}\frac{\partial}{\partial x_i}\end{equation*}    Besides, it is identified with the Schur function $s_{\tau}$ and $C_{\tau}(YX) = C_{\tau}(\sqrt{Y}X\sqrt{Y})$ for any Hermitian matrix $Y$.
It is well-known that, if $p=q+1$, then the hypergeometric series is convergent for $0 \leq ||X|| < 1$ ($|| \cdot ||$ is the norm given by the spectral radius) , if $p \leq q$, then it converges everywhere and else, it diverges.
\subsection{The modified Bessel function (Lebedev 1972)}
The modified Bessel function with index $\nu \in \re$ is given by the following series : \begin{equation*}{\it I}_{\nu}(z) = \sum_{k=0}^{\infty}\frac{1}{k!\Gamma(\nu+k+1)}\left(\frac{z}{2}\right)^{2k+\nu} ,\qquad z \in \mathbb{C}.\end{equation*} It can be represented through standard hypergeometric functions ${}_0F_1$ and ${}_1F_1$  : 
\begin{align*}{\it I}_{\nu}(z) &
= \frac{1}{\Gamma(\nu+1)}\left(\frac{z}{2}\right)^{\nu}{}_0F_1(\nu+1; z^2)
\end{align*}
{\bf Acknowledgments} I would like to thank {\it C. Donati-Martin} (my Ph. D. supervisor), {\it J. Faraut} and  {\it M. Yor} for  helpful comments and encouragements. 

\linespread{1.6}

\end{document}